\title {     
Sharp estimates for Brownian non-intersection probabilities
}
\author {Gregory F. Lawler\thanks{Duke University, Research
supported by the National Science Foundation}
\and Oded Schramm\thanks{Microsoft Research}
\and Wendelin Werner\thanks{Universit\'e Paris-Sud}
}
\documentclass[12pt,naturalnames]{article}
\usepackage{amsmath}
\usepackage{amsthm}
\usepackage{amsfonts}
\newif\ifhyper\IfFileExists{hyperref.sty}{\hypertrue}{\hyperfalse}
\ifhyper\usepackage{hyperref}\fi

\newif\ifdraft
\drafttrue
\numberwithin{equation}{section}

\newtheorem{theorem}{Theorem}
\numberwithin{theorem}{section}
\newtheorem{corollary}[theorem]{Corollary}
\newtheorem{lemma}[theorem]{Lemma}

\def \eps {\varepsilon}

\def \BM {Y}
\def \R {{\mathbb R}}

\def \expect {{\bf E}}
\def \P {{\bf P}}
\def \b {{\flat}}

\def \A {{\mathcal A}}
\def \C {{\mathcal C}}
\def \st {{\ : \ }}
\def \E {{\bf E}}

\def \D {{\mathcal D}}
\def \flat {{\delta}}

\def \Y {{\mathcal Y}}
\def \B {{\mathcal B}}

\usepackage{amsmath}
\usepackage {amsfonts}
\begin{document}
\maketitle

\begin{abstract}
This paper gives an accessible (but still technical)
self-contained proof to the fact
that the intersection probabilities for planar
Brownian motion are given in terms of the intersection
exponents, up to a bounded multiplicative error,
and some closely related results. 
While most of the results are already known, the proofs are
somewhat new, and the paper can serve as
a source for the estimates used in our paper 
\cite{LSWan} on the analyticity of the Brownian intersection 
exponents.
\end{abstract}

\section{Introduction}

In a recent series of papers \cite{LSW1,LSW2,LSW3,LSWan},
the authors rigorously derived the values for the intersection
exponents for planar Brownian motion. Among other things,
we prove in these papers  
that the Hausdorff dimension
of the outer boundary of a planar 
Brownian path is $4/3$ (see \cite {LSW4/3} for an overview).
This  paper
is complementary to these papers in that it proves
some estimates about the Brownian intersection probabilities
that do not depend on knowing their exact values.

The intersection exponents are defined as the 
asymptotic  rate
of decay of certain non-intersection probabilities.
The main results in this paper give estimates
for the actual probabilities in terms of these
asymptotic exponents.  For
example, it is easy to 
show by subbadditivity
that the probability that the paths of
two independent planar Brownian motions started
from uniform-independent points on the
unit circle will not intersect before hitting
the circle of radius $R$ is 
$R^{-\xi+o(1)}$ as $R\to\infty$
(this is  the definition of the
exponent $\xi = \xi(1,1)$).  We show that in fact
the probability is equal to $u(R)R^{-\xi}$,
where $u(R)\in[1/c,c]$ for some constant $c$ independent of $R \ge 1$.

This, and most other results 
proven here have been derived before  
by Lawler %%O-ch
 (see \cite{Lstrict}
and reference therein).  However, these earlier treatments
were a little complicated at times  (one reason is that
 they simultaneously
treated both the planar and three-dimensional cases) and
hence it seems worthwhile to have a self-contained account of
these results.   We will not make any use of our recent papers
\cite{LSW1,LSW2,LSWan,LSW3}; instead this paper can be
considered as a prerequisite to \cite {LSWan}.
The results presented here are used 
in \cite {LSWan} to prove analyticity of the mappings $ \lambda 
\mapsto \xi (k, \lambda)$ on $(0, \infty)$. 
``Up-to-constants'' estimates are
 also instrumental in relating the intersection exponent
to the Hausdorff dimension of exceptional sets of the
Brownian path, see, e.g., \cite{Lfrontier}.

We will concentrate on the intersection exponents
 $\xi(2,\lambda)$ which
are  relevant  for analyzing the outer boundary of
Brownian paths.  However, the proofs, with only minor changes,
adapt easily to other Brownian intersection exponents (see
Section
\ref{final}).

For all $r \ge 0$,
let $\C_r$ denote the circle %%O-ch
% define the circle $\C_r$
of radius $\exp (r)$ about zero.
Let $\BM^0,\BM^1,\BM^2$  be independent planar Brownian
motions starting at $0$.  Define for $j=0,1,2$,
and $r \in \R$,
\[  T^j_r = \inf\{t >0 \st  \BM^j_t \in \C_r  \}\]
and the paths
\[ \Y^j_r = Y^j [ T_0^j, T_r^j ] \]
(one could equivalently have taken Brownian motions
started uniformly on the unit circle up to their
hitting time of $\C_r$).
We define the random variable (depending on $\Y_r^1$
and $\Y^2_r$),
$$
Z_r = Z_r (\Y^1_r, \Y^2_r) 
:= 
\P [ 
 \Y_r^0 \cap ( \Y_r^1 \cup \Y_r^2 )= \emptyset \mid \Y^1_r, \Y^2_r].
$$
This is the probability, given $\Y_r^1$ and $\Y_r^2$,
 that another Brownian motion started
uniformly on the unit circle reaches $\C_r$ 
without intersecting the paths $\Y_r^1$ and $\Y^2_r$.
We define for all $\lambda>0$,
$$ 
a_r= a_r (\lambda)  = \E [ (Z_r )^\lambda ].
$$
Note that when $\lambda$ is an integer, then $a_r$ is the probability
that $\lambda$ independent copies of $\Y^0_r$ do not intersect $\Y_r^1 
\cup \Y_r^2$. It is straightforward
 to show that there exists a constant
$\xi$, usually denoted by $\xi (2 , \lambda)$, such that 
$$
\lim_{r \to \infty} (a_r)^{1/r}=  e^{- \xi} 
.$$
One of the main goals of the present paper is to present a 
short and self-contained proof of estimates for $a_r$
(and alternative closely related quantities) 
up to multiplicative constants. In particular:

\begin {theorem}
\label {main}
For every $\lambda_0 > 0$, there exist constants 
$0<c_1< c_2< \infty$ such that for every $0 <\lambda \le \lambda_0$
and every $r \ge2$,
$$
c_1 e^{-r \xi(2, \lambda) } 
\le
a_r (\lambda)
\le 
c_2 e^{-r \xi (2, \lambda)}.
$$
\end {theorem}

This theorem is a slight improvement over the estimate given
for $a_r$ in \cite{Lstrict}.  In that paper, it was
shown that for every $0 < \lambda_1 < \lambda_2 < \infty$, one can
find constants $c_1,c_2$ that work for all
$\lambda\in[\lambda_1 , \lambda_2]$.  The approach we give in this paper
gives the stronger result that the constants can be chosen
uniformly in $(0,\lambda_0]$.  An advantage of Theorem \ref{main}
is that the following is an easy corollary obtained by
fixing $r$ and letting $\lambda \rightarrow 0+$.

\begin{corollary}  There exist constants $0 < c_1 <c_2 < \infty$ 
such that
$$
     c_1 e^{-r \xi(2,0) } \leq  \P [ Z_r > 0 ]
    \leq c_2 e^{-r \xi(2,0)} ,
$$
where $\xi(2,0) :=\lim_{\lambda \rightarrow  0+} \xi(2,\lambda)$.
\end{corollary}
Note that $Z_r >0$ means that $\Y_r^1  \cup \Y_r^2$ does 
not disconnect 
$\C_0$ from $\C_r$.
This corollary  was derived in \cite{Lfrontier} for
the disconnection exponent $\xi_{0}$ defined by
$e^{-\xi_0}=\lim_{r\to\infty}\P [  Z_r > 0 ]^{1/r}$.
  However, a more complicated
argument was needed~\cite{Lstrict} to prove
$\xi_0 =\lim_{\lambda\rightarrow 0+} \xi(2,\lambda)$.  Using Theorem
\ref{main}, this is immediate.

Although we do not prove it in this paper, it can
actually be shown that  quantities
like $e^{r\xi(2,\lambda)} a_r$ approach %%O-ch
a limit (see the end of Section \ref{separatesec}).

Another goal of the present paper is to clarify 
and summarize the equivalence 
between the definitions of the exponents in terms of 
Brownian excursions, Brownian motions, extremal distance, %%O-ch length,
and discuss %%O-ch
the influence of the starting points, etc. 
In fact, we will  
first focus on another quantity $b_r$ defined in terms of 
Brownian excursions and extremal distance, 
show up-to-constants estimates for $b_r$ and then deduce 
the estimates for 
$a_r$.
 
\section {Preliminaries}

Before studying non-intersection probabilities, we first 
review a few easy facts concerning Brownian excursions
and
 extremal
distance.

Throughout the paper, for all $r<r'$, $\C_r$ will denote the
circle of radius $\exp (r)$ about $0$, and $\A (r,r')$ 
will denote the open annulus between $\C_r$ and $\C_{r'}$.
$\D(z,\delta)$ will denote the open disk of radius $\delta$ about $z$.
It will be sometimes more convenient to work in the
cylindrical metric. We will then implicitely use
the fact that for all $\eps>0$,
 when $\delta$ is sufficiently small, for
all  $z= e^{u} \in \C_r$, %%O-ch then 
$$
\D (z, \delta e^r (1- \eps)) \subset
\{ e^v \st | v - u| < \delta \} \subset \D( z, \delta e^r (1+\eps) ) .
$$

\subsection{Excursion measure and conformal invariance}

Let $Y$ be a Brownian motion starting
at the origin, let 
$T_r$ be its hitting time of the circle $\C_r$
and define 
\[   S_r = \sup\{t < T_r \st \BM_t \in \C_0\}.  \]
 The paths 
\begin{equation}  \label{dec12.1}
  B_t := Y_{t}  , \;\;\;\;
    S_r \leq t \leq T_r  , 
\end{equation}
are called ``Brownian upcrossings'' of the annulus $\A(0,r)$.
We will not care about the time-parameterization
of the upcrossings; in particular,
it does not matter if the `starting-time' of the upcrossings 
is called $0$ or $S_r$.

This probability measure
on  Brownian upcrossings is very closely related to 
the Brownian excursion measure that we used 
in the papers \cite {LW1,LW2,LSW1,LSW2}.  
The excursion measure on
the annulus $\A(0,r)$ is the upcrossing
probability
normalized so that
 the total mass is $2 \pi r^{-1} $.

 We now briefly 
recall some of the properties of these measures.
First, there are various equivalent ways of defining them.
Define
the excursion measure on $\A(0,r)$ starting at $z \in \C_0$
by \[\mu_{z,r} = \epsilon^{-1}
\lim_{\epsilon \rightarrow 0} \mu_{z,r,\epsilon} \] where
$\mu_{z,r,\epsilon}$ is the measure obtained
from starting a Brownian motion at
$(1 + \epsilon)z$,  killing it upon leaving $\A(0,r)$, and
restricting to those paths that exit $\A(0,r)$ at $\C_{r}$.
Then the excursion measure on $\A(0,r)$ is given by
\begin{equation}  \label{dec18.1}
\int_0^{2 \pi} \mu_{\exp (i \theta),r} \; d\theta .
\end{equation}

Yet another equivalent way to define the probability measure on 
upcrossings is to identify upcrossings
 with the process  $R_t= \exp (U^1_t+ i U^2_t)$ 
where $U^1$ is a three-dimensional Bessel process started at $0$, 
and $U^2$ an independent Brownian motion started uniformly on $[0, 2
\pi]$,
stopped at the first time it hits the circle $\C_r$ (i.e., %%O-ch
 at the first time
$U^1$ hits $r$) (see e.g. \cite {RY} for definition
and properties of Bessel processes).

When $r<r'$, define the excursion measure and the upcrossing 
probability on $\A(r,r')$ as the measure (or law) of $e^r$ %%O-ch
 times 
a Brownian upcrossing in $\A(0, r'-r)$.
It is easy to see (for instance using the definition of
the upcrossings in terms of Bessel processes) 
that if $B$ is a Brownian upcrossing of $\A(r,r')$,
then the time-reversal of $1/B$ is a Brownian upcrossing of 
$\A(-r', -r)$.

One can in fact define excursion measures in any
open planar domain. 
In the papers \cite {LW1,LW2,LSW1,LSW2} we used Brownian excursion
measures in simply connected planar domains. Just as in \cite {LW1},
in the present paper,
 we will need to use this measure only in
some particular simply connected
 domains.   Suppose $O$ is a simply connected subset
of $\A(r,r')$,
and define  $\partial_1 := \partial O
\cap \C_r$ and $\partial_2 := \partial O \cap \C_{r'}$.
Let $\Phi$ denote a conformal map from $O$ onto the unit disk.
We say that $O$ is a path domain 
in $\A (r,r')$ if $\Phi (\partial_1)$ and 
$\Phi (\partial_2)$ are two arcs of positive length.
We call $\partial_3$ and $\partial_4$ 
the two other parts  of $\partial O$ (possibly viewed as
sets of  prime ends).
The excursion measure in $O$ can be defined
as the excursion measure in $\A(r,r')$
restricted to those upcrossings that stay in $O$. 

An important
property of the excursion measure is its conformal invariance:  if
$F$ is a conformal transformation taking a path domain $O$
 to another path domain $O'$ in such a way that $F( \partial_1)
= \partial_1'$ and $F(\partial_2) = \partial_2'$ (with obvious 
notation)
then the image of the excursion measure on $O$
by $F$ is the excursion measure on $O'$. 
See for instance \cite {LW1,LW2} for a proof of this fact.

\subsection {Extremal distance and excursions}
\label {exsec}

For any path domain $O$, there exists a unique 
positive real $L$ such that $O$ can be mapped 
conformally onto the half-annulus $O_L'=
\{ \exp (u + i \theta) 
\st 0 < u < L \hbox { and } 0 < \theta < \pi \} $
in such a way that $\partial_1$
and $\partial_2$ are mapped onto the semi-circles
(or equivalently, such that $O$ can be mapped conformally onto
the rectangle ${\cal R}_L := (0, L) \times 
(0, \pi)$ in such a way that 
$\partial_1$ and $\partial_2$ are mapped
 onto the vertical sides of ${\cal R}_L$).
We call $L = L(O)$ the 
$\pi$-extremal
distance between $\partial_1$ and $\partial_2$ in $O$.
This is 
$\pi$ times the extremal distance %%O-ch length
as
in \cite{A,Pom}.

The excursion measure can also be defined on 
the rectangle  
${\cal R}_L$  by taking image of the excursion measure in $O_L'$ 
under the logarithmic map.
Alternatively, it can 
directly be defined as $\pi / \eps$ times the limit when $\eps \to 0$
of 
the law of Brownian paths started uniformly on the 
 segment $[\eps, \eps + i \pi ]$ and 
restricted to the event where they
 exit the rectangle through $[L, L+ i \pi]$.

Since the excursion
measure is invariant under conformal transformations,
its total mass depends only on $L$.
By considering directly excursions in the rectangle,
it is easy to check that
 there exists a constant $c$ such
that for all $L \geq 1$, the total mass of the excursion
measure in ${\cal R}_L$ is %%O-ch
 in $[c^{-1} e^{-L},
c e^{-L}].$  
In other words, 
up to multiplicative constants, $e^{-L(O)}$ measures the 
total mass
of the excursion measure in $O$.

Extremal distance in a simply connected domain $O$ can be defined
in a %%O-ch
 more general context. For instance (see, e.g., %%O-ch
\cite {A}), suppose
that $V_1$ and $V_2$ are arcs on the boundary of $O$, and let 
$\Gamma$ denote the set of (smooth) paths that disconnect $V_1$
from $V_2$ in $O$. For any piecewise smooth metric 
$\rho$ in $O$, define
the $\rho$-area  %%O-ch
$ {\cal A}_\rho (O) :=  \int_O \rho(x+iy)^2 \: dx \: dy $ of $O$ and
the length of smooth curves $\gamma$,
$       \ell_\rho (\gamma):= 
   \int_\gamma \rho(z) \: d|z|$.
Then, define
$$
L ( O ; V_1, V_2) : =  \pi \inf_{\rho} {\cal A}_\rho (O)
$$
where the infimum is taken over the set of piecewise smooth
metrics such that 
for all $\gamma \in \Gamma$, $\ell_\rho (\gamma) \ge 1$.
It is straightforward to check that this definition generalizes the
previous one (it is also invariant under 
conformal transformations, and if $V_1$ and $V_2$
are the vertical sides of $O={\cal R}_L$, the infimum is 
obtained for a constant $\rho = 1/\pi$).

Using rectangles, it is easy to see that this 
definition is equivalent to the more usual definition 
(see \cite {A}) of extremal distance in terms of 
the family of curves connecting $V_1$
to $V_2$ in $O$ (i.e. $L(O; V_1, V_2)$ is the
 maximum over all metrics 
$\rho$ with ${\cal A}_\rho (O) = 1$ of the square of 
the $\rho$-distance between $V_1$ and $V_2$ in $O$).

It is straightforward to see that $L(O;V_1, V_2)$
 satisfies monotonicity
relations: if $O' \subset
O, \partial_1' \subset \partial_1,$ and $\partial_2' \subset
\partial_2$, then ${L}(O';\partial_1',\partial_2') \geq
L(O;\partial_1,\partial_2)$;
and if $C$ is a  simple curve in $O$ connecting $\partial_3$ and
$\partial_4$, $O'$ is the connected component of
$O \setminus C$ whose boundary contains
 $\partial_1$, and
 $O^*$ is the  component of
$O \setminus C$ whose boundary contains $\partial_2$, then
$
L(O;\partial_1,\partial_2) \geq L(O';\partial_1,C)
+ L(O^*;C,\partial_2)$.

\subsection{A few simple lemmas}

We will need a few simple  technical facts about extremal
distance. It will be more convenient here
to work with the cylindrical metric.
 Let $\tilde O$ be a path domain on $\A(0,r)$ 
(with $\tilde \partial_1 , \ldots, \tilde \partial_4$ being the
four parts of $\partial \tilde O$).
Throughout this section, we will use a simply connected 
set $O$ such that $\exp (O) = \tilde O$.
We define $\partial_1 , \ldots,\partial_4$
the parts of $  O$
corresponding to $\tilde \partial_1, \ldots, \tilde 
\partial_4$, and we will suppose that $\partial_3$ is
`below' $\partial_4$ 
(i.e. that $z_1:= \partial_3 \cap \{ \Re (z) =0 \}$
lies below $z_2 : = \partial_4 \cap \{ \Re (z) = 0 \}$).
Note that 
$ \partial_3 \cap  \partial_4 = \emptyset$
(while it was possible that 
$\tilde \partial_3 \cap \tilde \partial_4
\not= \emptyset$).
The following lemmas will be formulated in terms
of $O$, and applied later to $\tilde O = \exp (O)$.
We will not bother to choose optimal constants
as only their existence  will be needed.

\begin{lemma}  \label{dec17.lemma3}
Suppose that 
for some $\delta<r-1$,
$
    \D(z_1,4 \delta) \cap \partial_4 = \emptyset$ and $
 \D(z_2,4 \delta) \cap \partial_3 = \emptyset .$
Then,
\[  L(O \setminus [\D(z_1, \delta)
  \cup \D(z_2, \delta)]; \partial_1,\partial_2) \leq
   L(O;\partial_1,\partial_2)  +  6 \pi^2 . \] 
\end{lemma}

{\bf Proof.}  Let $O'$ be the domain 
$O \setminus [\D(z_1, \delta)
  \cup \D(z_2, \delta)]$ and write $ \partial_1'
= [z_1 + i \delta, z_2 - i \delta],
\partial_2' = \partial_2, \partial_3', \partial_4'$
for the corresponding boundaries.  Let $\rho$ be the extremal
metric for finding the length of the collection $\Gamma$
of curves  in
$O$ connecting $\partial_3$ and $ \partial_4$ 
(note that $\rho$ is 
the conformal 
image of a multiple of the Euclidean metric in the rectangle, and
therefore $\rho$ is smooth)
so that
\[ {\cal A}_{ \rho} ( O) =  
    \pi^{-1} L(O;\partial_1,\partial_2). \]
If we let $\Gamma'$ be the collection of curves in $ O'$ connecting
$ \partial_3'$ and $ \partial_4'$, 
 and 
\[      \rho' = \max\{ \rho,
 \delta^{-1} [1_{\D(z_1,2 \delta)}
               + 1_{\D(z_2,2 \delta)} ]\} \]
in $ O'$, then every curve in $\Gamma'$ has length at least one in the
metric $ 
\rho'$.  Hence
\begin{eqnarray*}   L(O';\partial_1',\partial_2')
  &  \leq & \pi  {\cal A}_{\rho'} (O') \\
   &  \leq & \pi [{\cal A}_{\rho} (O) + 2 (4\pi - \pi)] \\ %%O-ch
   & = & L(O;\partial_1,\partial_2) + 6 \pi^2 . \qed %%O-ch
\end{eqnarray*}

\begin {lemma}
\label {jan5}
For all $\delta>0$,
there exists $c(\delta )$ such that
if $V 
\subset \partial_1$ is a segment of length at least 
$\delta$, if ${\rm dist}(V, \partial_3 \cup \partial_4) > 
 \delta$ and if the $\delta$-neighborhood of $V$ disconnects
$\partial_3$ from $\partial_4$ in $O \cap ((0, \delta) \times \R)$, then 
$$
L ( O ; V, \partial_2 ) \le L (O; \partial_1, \partial_2) + c(\delta).
$$
\end {lemma}

{\bf Proof.}
Let $ \rho$ denote the extremal metric in $O$ associated to $L(
O; \partial_1, \partial_2)$ (i.e., %%O-ch
 any path from $ \partial_3$
to $ \partial_4$ in $O$ has $ 
\rho$-length at least one, and ${\cal A}_{ \rho}
(O)$ is minimal), and define 
$$ \rho' =
\max \{ \rho , \delta^{-1} 1_{(0,  \delta ) \times \R } ].
$$
Any path disconnecting $V$ from $ 
\partial_2$ has $\rho'$-length
at least one, so that 
$L (O; V, \partial_2) \le \pi \A_{ \rho'} ( O) $
and the lemma follows.\qed

\begin {lemma}
 \label{dec17.lemma2}
\label {lastone}
Suppose  that $1<s<r-1$,  
and that for some  small $\delta$,  
$\partial_3' :=\partial_3 \cap ((s- \delta, s+ \delta) \times \R)$ and
$\partial_4':= \partial_4 \cap ((s- \delta , s+ \delta) \times \R)$ 
are both of diameter smaller than $\delta^{1/6}$ and at distance 
at least $\delta^{1/7}$ from each other.
Let $V$ denote the segment in $O \cap \{ Re (z)=s \}$
that disconnects $\partial_1$ from $\partial_2$ (it is unique 
because of the previous conditions).
Then,  for some $C (\delta)$, 
$$
L (O ; \partial_1, \partial_2) 
\le L ( O \cap ((0,s) \times \R) ; \partial_1 , V) 
+ L (O \cap ((s,r) \times \R ) ; V , \partial_2 ) 
+ C(\delta).
$$
\end {lemma}

{\bf Proof.}
Let $O_1$ and $O_2$ denote the sets
$O \cap ((0,s) \times \R)$ and $O \cap ((s,r) \times \R)$.
Let $\rho_1$ (resp., $\rho_2$) denote the extremal metric 
in $ O_1$ associated to
$ L(O_1 ; \partial_1,V)$ (resp., in $
O_2$ 
associated to $ L(O_2;
V,\partial_2)$). Let ${\cal V} 
=  O \cap (( s-\delta, s+ \delta)\times \R)$. Note that 
(since $\exp O = \tilde O$) the euclidean area of ${\cal V}$
is at most $4 \pi \delta$.
Define
$$
\rho  = 
\max (  \rho_1  ,  \rho_2  , (1/\delta) 1_{
{\cal V}} ).$$
It is easy to check that any path joining 
$\partial_3$ to $ \partial_4$
in $O$ has $\rho$-length at least $1$ 
(either, it stays in one of the three sets
$ O_1$, 
 $ O_2$ or ${\cal V}$, or it contains a
path joining $\{ \Re (z) = s \}$
to $\{ |\Re (z) - s | = \delta  \}$).
Therefore, 
$$
 L ( O; \partial_1, \partial_2 ) \le \pi {\cal A}_{\rho} (O)
\le  L(O_1 ; \partial_1,V) + L(O_2;
    V,\partial_2) + C(\delta).
\qed
$$

\subsection {Extending excursions}

Let $0<r<r'$.
A consequence of the strong Markov property of planar 
Brownian motion and of the second definition of the 
Brownian excursion measure is that if $B$
is a Brownian upcrossing of $\A(0, r)$ defined
under the excursion measure, and if one starts from its 
endpoint (on $\C_r$) another independent planar
Brownian motion killed at its first 
hitting of $\C_{r'}$, restricted to the event that it does
not intersect $\C_0$ (note that 
this is an event of probability $r/r'$),
then the concatenation of the upcrossing with the Brownian path
is exactly defined under the Brownian excursion measure in $\A(0,r')$.
 
In particular, this shows that if $B$ is an 
Brownian upcrossing of $\A(0,r')$
(defined under the probability measure on upcrossings), then
it can be decomposed into two parts: A Brownian upcrossing 
of $\A (0,r)$ and a Brownian motion started from the end-point of the 
first part, conditioned to hit $\C_{r'}$ before $\C_0$.

This can also be formulated easily in terms of the definition of 
Brownian upcrossings using three-dimensional Bessel processes.
In particular, 
it shows that it is possible to define on the same probability 
space a process $(R_t, t \ge 0)$ started uniformly on the unit
circle, such that for each $r>0$, the process $R$ stopped at its 
hitting time $T_r$ of the circle $\C_r$ is a Brownian upcrossing  
of $\A(0,r)$. We will use the $\sigma$-field  ${\cal F}_r$ 
generated by $(R_t, t \le T_r)$ in Section \ref {separatesec}.

Another simple consequence of the strong Markov property of 
planar Brownian motion is the fact that 
conditionaly on $Y(S_r)$ (which has 
uniform law on $\C_0$), the Brownian
upcrossing $Y [S_r , T_r]$ is independent from the initial 
part $Y [0, S_r]$.
Consider now the event 
$H = H_r$ that $Y[T_0,S_r]$ 
does not contain a closed
loop about zero contained entirely in the annulus
$\A(-1,0)$. This event is independent of the 
upcrossing $Y[S_r, T_r]$ so that on this event, 
the measure on upcrossings
is the same as the
upcrossing probability or the excursion measure except that it
has a slightly different normalization constant i.e., %%O-ch
its total mass $m_r$ is the probability of $H_r$.
We claim there is a constant $c$ such that
 $ c^{-1} r^{-1}
\leq m_r \leq c r^{-1}$.  The lower bound can for instance be
derived by considering the event $\{Y[T_0,
T_r] \cap \A(-1,0) \subset \D(Y(T_0), \delta)\}$
for some fixed $\delta <1/4$.
For the upper bound, let $k$ denote the total number
of times the Brownian motion goes from $\C_0$ to
$\C_{-1}$ before time $T_r$.  Every time the
path goes from $\C_0$ to $\C_{-1}$ there is a positive
probability, say $\rho$ of forming a closed loop in
$\A(-1,0)$.  From this and the strong Markov property,
we get $\P(H_r \cap \{k = l\}) \leq (1-\rho)^l r^{-1}$,
and summing over $l$ gives the upper bound.

We note that we have just proved that
for all $\delta <1/4$, there is a $c'=
c'(\delta)$ such that conditioned
on the event $H_r$, the probability that $Y[T_0,
T_r] \cap \A(-1,0) \subset \D(Y(T_0), \delta)$
is at least $c'$.

\section{Lower bound}  

From this point on, we fix a $\lambda_0$ and consider
$\lambda \in (0,\lambda_0]$.  Constants are allowed
to depend on $\lambda_0$ but not on $\lambda$.

Suppose  that $B^1$ and $B^2$ are  two 
independent Brownian upcrossings of the annulus $\A (0,r)$ defined
using the Brownian motions $Y^1$ and $Y^2$.
Let $O^1$ and $O^2$ be the components of  $\A (0,r) 
\setminus (B^1 \cup B^2)$ which are at zero distance from $\C_r$.
We choose $O^1$ in such a way that it has the positively 
oriented arc on $\C_r$
from the endpoint of $B^1$ to the endpoint
of $B^2$ as part of its boundary.
For $j=1,2$, let  $L^j_r=L (O^j)$ be 
 $\pi$ times the extremal 
distance between 
$\C_0 \cap \partial O^j$ and $\C_r \cap \partial O^j$ in
 $O^j$. 
Note that  $\C_0 \cap \partial O^j$ is a.s.\ either empty
or an arc, and $ \C_r \cap \partial O^j$ is a.s. an arc (note that in
this case $O^j$ is a.s. a path domain).
When $\partial O^j\cap\C_0=\emptyset$, set ${ L}^j_r:=\infty$.
Let ${ L}_r := \min \{ L^1_r, L^2_r \}$, and 
 let $O:=O^j$ when $L_r = L^j_r < \infty$.
Define
$$
b_r = b_r (\lambda)
= 
r^{-2} \E [ \exp (- \lambda {L}_r ) ].
$$

The goal of the next two sections is to define the intersection
exponent $\xi (2, \lambda)$ in terms of $b_r$, and to 
prove the following estimates for $b_r$. 

\begin {theorem}
\label {main2}
For any $\lambda >0$, there exists $\xi (2, \lambda) \in (0, \infty)$
such that $e^{- \xi (2, \lambda)} = \lim_{r \to \infty}
b_r^{1/r}$.
Furthermore, for any $\lambda_0>0$,
 there exist constants $c_1$ and $c_2$ such that 
for all $\lambda \in (0, \lambda_0]$, and for all $r \ge 0$,
$$
c_1 e^{- r \xi (2, \lambda)} \le b_r (\lambda) \le
c_2 e^{-r \xi (2, \lambda)}.
$$
\end {theorem}

In the present section, we will derive the lower 
bound
and the next section will be devoted to the (harder) upper
bound, and we will relate $a_r$ to $b_r$ in the subsequent 
section.
Note that $b_r$ is decreasing in $r$
because of the monotonicity properties of extremal distance.

For any positive integer $n$,
let $E_n$ denote the event that neither $\BM^1
[T_0^1,T_n^1]$ nor $\BM^2[T_0^2,T_n^2]$ 
hit the circle $\C_{-1}$.
Note that $\P (E_n) = 1/ (n+1)^2$ and that $E_n$ is independent
from $\BM^1 [ S_n^1 , T_n^1]$ and $\BM^2 [S_n^2, T_n^2]$.
Hence
$$ 
E [ e^{-\lambda L_n} 1_{E_n} ] =
n^2 b_n / (n+1)^2.$$
We call $b_n^\#$ this quantity.

\begin{lemma}  \label{lemma.dec9.1}
There exists a constant $c$ such
that for all $n,m \geq 1$,
\[   b_{m+n+1} \leq c b_m b_n. \]
\end{lemma}

{\bf Proof.}
First consider the event $H_n^1 \cap H_n^2$ that neither
$Y^1[T_0^1,T_n^1]$ nor $
Y^2[T_0^2,T_n^2]$ contains a closed loop in
$\A(-1,0)$ that surrounds $\C_{-1}$.  
The previous considerations show that 
$L_n$ is independent from $H_n^1 \cap H_n^2$
so that 
\begin{equation}  \label{dec7.1}
   b_n^* := \E[e^{-\lambda L_n} 1_{H_n^1 \cap H_n^2} ]  \leq 
\P [ H_n^1 \cap H_n^2 ] \E [ e^{-\lambda L_n} ]
\le c b_n .
\end {equation}
Once we have this, to get the lemma we split the
upcrossings into the pieces up to $T_m^j$ and from
$T_{m+1}^j$ to $T_{m+n+1}^j$.  Monotonicity of
extremal distance gives
\[    b_{m+n+1}^\# \leq c b_m^\# b_n^*, \]
from which the lemma follows.  \qed

\vspace{4ex}

Using this lemma, we can now define $\xi(2,\lambda)$
by
$  e^{-\xi(2,\lambda)} = 
     \lim_{n \rightarrow \infty}
          b_n^{1/n} $ 
and get $b_n \geq c e^{-n \xi(2,\lambda)}$
for some $c$, which gives the lower bound in
Theorem \ref{main2} for integer $n$'s.
    By considering Brownian motions restricted
to stay in the upper or lower half-plane we get the
crude estimate
$   \xi(2,\lambda) \leq 2 + \lambda \leq 2 + \lambda_0$.
We will use this fact implicitely in our estimates
when we write $e^{-\xi(2,\lambda)} \geq c$.  This is obvious, but
it is important that the constant can be chosen uniformly for
$0 < \lambda \leq \lambda_0$.  In this case $c=e^{-(2 + \lambda_0)}$
suffices.
In particular, since $b_r$ is decreasing in $r$, it suffices to 
prove the theorem for integer values of $r$.

\vspace{2ex}

In Section \ref{nicesec}
 we will need the following lemma.  Since 
the proof is very similar to that of (\ref{dec7.1}) we
include it here.  If $\epsilon > 0$, let $E_{n,\epsilon}$ be the
event that neither Brownian motion hits $\C_{-1 + \epsilon}$ 
before reaching $\C_n$.

\begin{lemma}  \label{lemma.dec9.2}  There is a constant 
$c$ such that for every $\epsilon \in (0,1/4)$,
\[  \E\: [ \: e^{-\lambda L_n} \: (1_{E_n} - 1_{E_{n,\epsilon}})
   \: ]
   \leq c \epsilon b_n^\#. \]
\end{lemma}

{\bf Proof.}  First note that $E_n \setminus E_{n, \epsilon}$ 
is independent of $L_n$ so that 
the left-hand side is equal to $\P [ E_n \setminus E_{n, \eps}]
\E[e^{-\lambda 
L_n}]$. Moreover, 
$\P [E_n] = 1/(n+1)^2$ 
and 
$\P [E_{n, \epsilon}]= (1-\eps)^2/(n+1-\eps)^2 
\ge ( 1 - c\eps) \P [E_n]$.
\qed

\section {The upper bound} 

Our goal in this section
is to derive the upper bound in Theorem \ref {main2}.
It suffices to give an upper bound
for  $\tilde{b}_n := n^{-2} \E[\exp(-\lambda L_n^1)]$ since
$\tilde b_n \le b_n \leq 2 \tilde{b}_n$. 
The basic strategy is to find a sequence
$(b_n^\delta)_{n \ge 1}$
such that:
\begin {itemize}
\item
For all $n \ge 1$, 
$b_n^\delta \leq \tilde{b}_n$.
\item There is a $c_1$ such that for all $n,m \ge 1$,
\begin{equation}  \label{dec9.1}
b_{n+m+2}^\delta \geq
c_1 b_n^\delta b_m^\delta.
\end{equation}
\item For all $n\ge 1$,
\begin{equation} \label{dec7.3}
   \#\{j \in \{1,\ldots,n\} : b_{j}^\delta \geq
    \tilde{b}_j /2 \} \geq 3n/4  . 
\end{equation}
\end{itemize}
Suppose we find such a sequence $b_n^\delta$.
It is then easy to check that
$\lim_{n \to \infty}
 (b_n^\delta)^{1/n} = \lim_{n \to \infty} (b_n)^{1/n} 
  =  e^{-\xi(2,\lambda)}$, and using (\ref{dec9.1}), that there
is a constant $c_3$ such that   for all $n \ge 1$,
\[    b_n^\delta \leq c_3  e^{-n \xi(2,\lambda)} . \]
Also, (\ref{dec7.3})	implies that for each $n \geq 2$ there is
a $j  \in \{1,\ldots,n-1\}$ such that
\[    b_j^\delta \geq  b_j/2, \;\;\; b_{n-j}^\delta \geq
      b_{n-j}/2 , \]
and hence
\[    \tilde{b}_{n+2}\leq b_{n+1} \leq c b_{j} b_{n-j} \leq
       4   c  b_j^\delta b_{n-j}^\delta 
      \leq 4 c c_1^{-1}   b_{n+2}^\delta
 \leq 4  c_3 c c_1^{-1}         e^{-\xi(2,\lambda)(n+2)} . \]
This establishes the upper bound.

\subsection{Nice configurations}  \label{nicesec}

Throughout this section, we will use Brownian upcrossings
$B^1$ and $B^2$ of annuli $\A (r, r')$. For convenience, 
we use the convention that $B^j$ is started at time zero on
$\C_r$, and that $T_{r''}^j$ denotes the first time at which $B^j$
hits $\C_{r''}$.

We  define a class of ``nice'' configurations
for pairs of Brownian upcrossings $B^1,B^2$ of
$\A(r,r')$ for $r'-r>1$.  More precisely, we say that
the configuration is $\delta$-nice at the beginning
if: 
\begin{itemize}
\item  $L^1 < \infty$;
\item $d(B^1 (0), B^2 (0))
>  \delta^{1/8}e^r$.
\item
For all $\eta<\delta$,
$B^j [0, T^j_{r+\eta^{1/2}}] \subset \D (B^j(0), \eta^{1/4}
e^r)$ for $j=1,2$.
\item 
For all $\eta<\delta$,
$B^j [T^j_{r + \eta^{1/2}}, T_{r+1}] \cap
\A (r, r+ 4 \eta ) = \emptyset$ for $j=1,2$.
\item $B^{j}[T_{r+1}^j,T_{r'}^j] \cap \A(r,r+ 4 \delta) 
   = \emptyset$ for $j=1,2. $
\end{itemize}
Here we write $L^1 = L^1 (r,r')$ for the appropriate
$\pi$-extremal distance.
Note that (and this is the reason for which we
introduce conditions with $\eta<\delta$) if 
a domain is $\delta$-nice at the beginning, then it is 
$\delta'$-nice at the beginning
for any $\delta' < \delta$.

Note also that the second, third and fourth conditions are only
on $B^1[0,T_{r+1}^1]$ and $B^2[0,T_{r+1}^2]$.
If we use $U(\delta)$ to denote the event that all these three
conditions
 hold, then, as the law of
$B^j [ 0, T_{r+1}^j]$ is that of a Brownian upcrossing
of $\A (r, r+1)$, we get easily that 
\begin {equation}
\label {Usmall}
 \P[U(\delta)] \rightarrow 1 
\end {equation}
as
$\delta \rightarrow 0+$, uniformly in $r'> r+1$.
In particular, almost surely, the configuration of
a pair of Brownian upcrossings
is $\delta$-nice at the beginning for sufficiently small 
$\delta$.

Analogously, we can define the notion of ``$\delta$-nice
at the end'' and we say that the configuration is $\delta$-nice
if it is $\delta$-nice at the beginning and at the end. 

Suppose now that $B^1$ and $B^2$ are two independent Brownian
upcrossings of $\A (0,n)$. 
Note that when the configuration is $\delta$-nice, then one can 
find a subarc of length at least $\delta$ on 
$\C_0 \cap \partial O_n^1$ 
that satisfies the conditions of Lemma \ref {jan5}.
Also, $O_n^1$ satisfies the conditions of Lemma \ref {dec17.lemma3}.
We shall use this later on.

Let
\[  b_n^\delta = n^{-2} \E[e^{-\lambda L^1_n} 1_{\delta-\mbox{nice}} ] .
\]

\begin{lemma}  \label{lemma.dec9.4}
For every $\epsilon > 0$, there is a $\delta_0 > 0$ such
that for all $\delta \in (0,\delta_0)$,
\[             \tilde{b}_n - b_n^\delta  \leq \epsilon
             b_{n-2} . \]
\end{lemma}

{\bf Proof.} Let $V = V_{n,\delta}$ be the event that the
configuration is not $\delta$-nice at the beginning, and
let $U = U_{n,\delta}$ be the $U(\delta)$ as above.  By symmetry and the
time-reversal property of upcrossings, 
it suffices to show that for all $\delta$ sufficiently small,
\[      n^{-2} \E[e^{-\lambda L^1_n} 1_{V}] \leq \frac{\epsilon}{2}
b_{n-2} . \]
Note that $V \cap \{L^1_n < \infty \} \subset U^c \cup V_1$ where
\[   V_1 = \bigcup_{j=1}^2 \{B[^j T_1^j,T_n^j] \cap \A(0, 4 \delta)
   \neq \emptyset \} . \]
The strong Markov property, decompositions
of Brownian upcrossings and monotonicity of
extremal distance, combined with (\ref {Usmall})
imply that
\[      n^{-2} \E[e^{-\lambda L^1_n} 1_{U^c}] \leq c \P(U^c) b_{n-2} .
\]
On the other hand,
Lemma \ref{lemma.dec9.2} establishes that
\[  n^{-2}  \E[e^{-\lambda L^1_n} 1_{V_1}] \leq c \delta b_{n-1}
\leq c \delta b_{n-2} . \qed
  \]

\begin{corollary}  For all $\delta$ sufficiently small,
(\ref{dec7.3}) holds.
\end{corollary}

{\bf Proof.}
First, we claim that for all $n$ sufficiently large
\begin{equation}  \label{dec9.4}
 \#\{j \in \{1,\ldots,n\}: b_{j+2} \geq  c      b_j \} \geq .9 n , 
\end{equation}
where $c = e^{-80 (2 + \lambda_0)}
$.
To see this, assume not.  Then, 
for infinitely many $n$'s,
 there exists at least  $.05 n$ exceptional
even values or at least $.05 n$ exceptional odd values $j$ 
in $\{1, \cdots, n\}$ scuh that $b_{j+2} \le c b_j$ in which case
\[   b_{n+2} \leq e^{-80 (2 + \lambda_0)(.05n)} \leq e^{-2(2 +
\lambda_0)n}
 \le   e^{-2 \xi(2,\lambda) n }   \]
and this contradicts the lower bound on $b_{n+2}$.  By changing
the value of $c$, we can conclude that 
(\ref{dec9.4}) in fact holds for all $n \ge 1$.  Hence,
Lemma \ref{lemma.dec9.4} (for $\eps = c/4$)
implies that for 
all $\delta$ sufficiently small, at least  $90\%$ of the integers
$j$ in $\{1, \ldots, n \}$,
\[ \tilde{b}_j - b_j^\delta  \leq  c b_{j-2} / 4
\le b_j / 4  \le 
             {\tilde b}_{j}/2   \]
 so that $b_j^\delta \geq
\tilde{b}_j/2  $. \qed

\subsection{Pasting}  \label{sepdefsec}

The goal is now to paste together nice configurations 
in order to get a lower bound for $b_{n+m+2}^\delta$ in terms
of $b_n^\delta$ and $b_m^\delta$.
In order to do this, we will define ``very nice configurations''.

From now on, we fix a small value of $\delta$ such that 
(\ref {dec7.3}) holds.
We say that
a configuration of pairs of upcrossings $(B^1, B^2)$
of $\A(r,r')$ is ``very nice at the end'' if
\begin{itemize}
\item $L^1 < \infty$;
\item $B^j(T_{r'-(1/3)}^j, T_{r'}^j) \subset
   \A(r'- (1/2), r'), \;\;\; j=1,2. $
\item $B^1 \cap \A(r'-\frac{1}{5}, r') \subset \{z:
            -\frac{1}{10} \leq  \arg(z) \leq
    \frac{1}{10} \} ; $
\item $B^2 \cap \A(r'-\frac{1}{5}, r') \subset \{z:
           - \frac{1}{10} \leq |\arg(z) - \pi| \leq
        \frac{1}{10}\} .$
\item $|\arg (B^1 (T_{r'}^1)) | \le 1/20$, $|\arg (B^2 (T_{r'}^2)) - \pi
| 
\le 1/20 $.
 \end{itemize}
Note that there is no $\delta$ in this definition. 
Let
\[   \beta_n^\delta = n^{-2}
 \E[\exp\{-\lambda L^1_{n}\} 1_{\delta-{\hbox {nice
 at the  beginning
     and very nice at the end}}}] . \]
However, by symmetry, the expectation is the same if we require
the configuration to be  $\delta$-nice at the end and
``very nice at the beginning.''
The goal is to paste together some configurations in
$\A (0, n+1)$ that are ``very nice 
 at the end'' with configurations in $\A (n+1, n+m+2)$ that are 
``very nice at the beginning.''

Suppose that $z^1$ and $z^2$
are on $\C_{n+1}$, and let 
 us now define $\beta_{n+1}^\delta (z^1,z^2)$ just as
$\beta_{n+1}^\delta$
except that the upcrossings are conditionned to end at $z^1$ and $z^2$.
In particular, since the law of the endpoints are uniform on $\C_{n+1}$,
$\beta_{n+1}^\delta$ is the mean of $\beta_{n+1}^\delta (z^1, z^2)$,
when $z^1$ and $z^2$ are integrated over $\C_{n+1} \times \C_{n+1}$.
Note that $\beta_{n+1}^{\delta} (z^1, z^2) = 0$ as soon as 
$(z^1, z^2) \notin 
Q:= \{e^{n+1+ i \theta } \st |\theta| < 1/20 \} \times
\{ e^{n+1+ i \theta} \st |\theta - \pi | < 1/20 \} $.

If $\alpha \in (0,\pi)$,
the probability that
a complex Brownian motion starting at $\epsilon
\in (0,1)$ reaches the unit circle without leaving the
wedge $\{z: |\arg(z)| \leq \alpha\}$ is at least
$\epsilon^{\pi/(2\alpha)}$ (this is easy
for $\alpha = \pi/2$ and can be established for
other $\alpha$ by considering the map $z \mapsto
z^{\pi/(2\alpha)}$).  
Such considerations show easily that if the configuration
of upcrossings of $\A (0,n)$ is 
$\delta$-nice, then with probability 
at least $c' \delta^c$, one can extend the upcrossings up to the 
circle $\C_{n+1}$ in such a way that the extensions
first remain in different wedges (and also leave an empty wedge
between them),  that all the wedges intersect 
$\A ( n-1, n)$ only inside the disks of radius $\delta$ around
the points $B^j (T_n^j)$    
and
such that the obtained configuration 
of upcrossings of $\A (0, n+1)$ is very nice at the end.
Furthermore, Lemmas
\ref {dec17.lemma3},
\ref {dec17.lemma2} and \ref {jan5} show that we can also impose that 
$e^{-L_{n+1}^1} \ge c' e^{-L_n^1} \delta^c$.
Finally, note that the weighted densities of the endpoints (on these 
configurations) 
on $\C_{n+1}$ are bounded away from zero
  on  $Q$.
Combining all this, we get that for any $(z^1, z^2) \in Q$,
\begin {equation}
\label {wed}
\beta_{n+1}^\delta (z^1, z^2)
 \ge 
c' \delta^c b_n^{\delta}
= c'' b_n^{\delta}
\end {equation}
(recall that $\delta$ is fixed).
Now we consider Brownian upcrossings $B^1$ and $B^2$
 of $\A (0, n+1+1+n')$
that are decomposed as follows:
A Brownian upcrossing of $\A (0, n+1)$, an intermediate
part and a final
Brownian upcrossing of $\A (n+1, n+1+1+n')$.
By restricting ourselves only to the cases where 
the first parts create a $\delta$-nice configuration at the 
beginning and are very nice at the end, where
the intermediate parts are of diameter smaller than $\delta e^{n+1} /10$
and where the final parts are very nice at the beginning and
$\delta$-nice at the end, using Lemmas 
\ref {dec17.lemma3},
\ref {dec17.lemma2} and \ref {jan5} again,
we get that 
$$
b_{n+n'+2}^\delta \ge c b_n^\delta b_{n'}^\delta
$$
for some $c(\delta)$ (we omit the details here).
This establishes (\ref{dec9.1})
and finishes the proof of
the upper bound.
\qed

\section {Non-intersection 
probabilities}
\label {otherdef}

We now show how the 
preceeding results (and in particular the strong
approximation for $b_n^\flat$) can be used to derive
Theorem \ref {main}, and ``up-to-constants estimates''
for other quantities closely related to $a_n$ and $b_n$. 

\medbreak

{\bf Proof of Theorem \ref {main}.}
Let $\B_r^j = Y^j [ S_r^j, T_r^j]$
denote the traces of the upcrossings.
For the upper bound, it suffices for example
to remark that 
\begin {eqnarray*}
Z_r (\Y^1_r, \Y^2_r)
 &\le& \P [ \B_r^0 \cap (\B_r^1 \cup \B_r^2)
= \emptyset \mid \B^1, \B^2]\\
&&  \times \P [ \BM^0 [T_0^0, S_r^0] \hbox { does not disconnect } \C_0
\hbox { from infinity}]
\\
&&
\times 1_{\BM^1[T_0^1,S^1_r] \hbox { and } \BM^2 [T_0^2, S^2_r] \hbox {
do 
not disconnect } \C_0 \hbox { from } \C_r}.  \\
\end {eqnarray*}
The first term is bounded above by  $c r e^{-L}$.
The second term is  
 bounded by a constant times $1/r$.
The last event is independent of $L$ and has probability
bounded by $c r^{-2}$. Therefore $\E[Z_r^{\lambda}] \leq
c b_r$.

For the lower bound, it suffices to use the lower bound for
$ b_n^\b$ 
and to realize the Brownian paths $\Y^j$'s  using a Brownian crossing
of the annulus together with initial parts 
 $Y^j [T_0^j, S_n^j]$ of small diameter.
\qed

\medbreak

Note that the estimates (3.1)
and (3.5) of \cite {LSWan} follow similarly.
Analogously, one can derive up-to-constants estimates if we prescribe
the
starting points of $\BM^1, \BM^2$ and/or of $\BM^0$ on the 
unit circle.
For instance, 
if we define
$$
\hat Z_n ( \Y^1_n, \Y^2_n ) 
= \sup_{z \in \C_0} 
\P [ \Y_n^0 \cap ( \Y_n^1 \cup \Y_n^2) 
= \emptyset \mid \BM^0 ( T_0^0) = z, \Y^1_n, \Y^2_n  ]$$
and
$$
\hat a_n 
= \sup_{z_1, z_2 \in \C_0 }
\E [(\hat Z_n )^{\lambda} \mid \BM^1 (T_0^1)= z_1, \BM^2 (T_0^2) =
z_2],$$
Then $a_n \leq \hat a_n$ and a simple application of
the strong Markov property  shows that $\hat a_n \leq c a_{n-1}$.
In particular,
$$
c_1' e^{-n \xi (2, \lambda)}
\le \hat a_n \le 
 c_2' e^{-n \xi (2, \lambda)},$$
for appropriately chosen $c_1',c_2'$.

\section {Separation lemma}  \label{separatesec}

In this section, we prove an important lemma that states that no matter 
how bad $O_n^1$ is, then there is a good chance 
(with respect to the normalized measure weighted by $\exp (-\lambda
L_{n+1}^1)$) 
that $O_{n+1}^1$ is very nice at the end as defined in Section 
\ref{sepdefsec}.
This lemma was the starting point for previous proofs of 
`up-to-constants' estimates, see \cite{Lstrict}.
While we do not need this lemma to establish the estimates
in this paper, we do use the 
 lemma in \cite {LSWan} to prove
analyticity of $\lambda \mapsto \xi (2, \lambda)$
(which was used to determine 
the disconnection exponents).  For this reason, we include
a proof here.

We use the notation of Section \ref{sepdefsec}. 
We suppose that the upcrossings $B^1, B^2$
 of $\A(0,r)$ are defined
in a compatible way in terms of Bessel processes
i.e., %%O-ch
 that both $B^1$ and $B^2$ are defined up to 
infinite time and that the upcrossings  $B^1 ( 0, T_r^1)$
and $B^2 (0, T_r^2)$ define the configuration
at radius $e^r$ ($O_r^1$ and $L_r^1$ are then defined in 
terms of these configurations). 
${\cal F}_r$ will denote the $\sigma$-field generated by these
two paths.
Recall that for all $r'>r$, conditionaly on ${\cal F}_r$,
the law of $B^j [T_r^j, T_{r'}^j]$
is that of a Brownian motion started from $B^j (T_r^j)$
conditioned to hit $\C_{r'}$ before $\C_0$.

Define
the event $\Delta(r,\delta)$  that the configuration
in $\A(0,r)$
in $\delta$-nice at the end, and the event $G_r$
that it is very nice at the end.

\begin {lemma}[Separation Lemma]
There exists  $c > 0$ such that for all 
$n \ge 1$, for all $\lambda \in (0, \lambda_0]$,
\begin {equation}
\label {sep}
\E [
1_{{G_{n+1}}} e^{- \lambda L_{n+1}^1 } \mid {\cal F}_n ]
\ge 
c 
\E [
 e^{-\lambda L_{n+1}^1}  \mid {\cal F}_n ].
\end {equation}
\end {lemma}

{\bf Proof.}  We start by noting that estimates for Brownian
motion in wedges show, just as for (\ref {wed}), 
that there exist $c,c'$ such that for any `stopping radius'
$\tau$ (i.e., %%O-ch
 stopping time for the filtration 
$({\cal F}_s)_{s \ge 0}$),
such that $\tau \in [n, n+1/4]$ almost surely,
\begin{equation}  \label{dec14.6}
 \E [
1_{{G_{n+1}}} e^{- \lambda L_{n+1}^1 } \mid {\cal F}_{\tau} ]
\ge 
c' \delta^{c} 
 e^{-\lambda L_{\tau}^1}  
 1_{\Delta(\tau,\delta)} 
\end{equation}
(because if the configuration is $\delta$-nice
at radius $\tau$, then one can extend it in such 
way 
that it is very nice at radius $n+1$).
Hence it suffices to find  $\delta_0, c''$ and such a  
stopping radius $\tau$ such that
\begin{equation}  \label{dec14.3}
  \E[ e^{-\lambda L_{ \tau}^1}
              1_{\Delta(\tau,\delta_0)} \mid {\cal F}_{n}]
  \geq c'' \E[e^{-\lambda L_{\tau}^1} \mid {\cal F}_n] . 
\end{equation}
For any positive integer $m$, let
\[  \tau_m = \inf\{s \geq 0: L_{n+s}^1 = \infty \mbox{ or }
          \Delta(n+s,2^{-m})   \} . \]
Note that if $L_n^1 < \infty$, then (up to 
a set of zero probability) $\tau_l = 0$ for all large enough $l$.

From the definition of $\delta$-nice
configurations, it is not difficult to see that
there exists $m_0$ and  $\rho > 0$ such that for all $m \geq m_0$,
\[  \P[\tau_m  \leq 2^{-m/20} \mid {\cal F}_n] \geq \rho . \]
By iterating, we see that 
\[  \P[\tau_m \geq  m^2 2^{-m/20} \mid {\cal F}_n] \leq
                      e^{-a m^2} , \]
for some positive constant $a$, and hence for all $m \geq m_0$,
\[ \E[e^{-\lambda L_{ \tau_m}^1} 1_{ \tau_m \geq  m^2 2^{-m/20}}
     \mid {\cal F}_n ] \leq   e^{-a m^2} e^{-\lambda 
L_n^1} . \]
On the other hand, using estimates in wedges again, we
see that for $m \ge m_0$, 
$$
\E [ e^{-\lambda L_{\tau_m}^1} \mid {\cal F}_n ]
\ge c 2^{-a'm} 1_{\Delta (n, 2^{-m-1})}  
 e^{-\lambda L_{n}^1} 
$$
so that
there is a summable sequence $\{h_m\}$ such that
\[ 
 1_{\Delta(n,2^{-(m+1)})}
\E[
e^{-\lambda L_{ \tau_m}^1} 1_{ \tau_m \geq   m^2 2^{-m/20}}
     \mid {\cal F}_n ]   \leq   h_m
    \E[
e^{-\lambda L_{ \tau_m}^1} 
     \mid {\cal F}_n ]  . \]
Similarly (starting at radius
$n+\tau_{m+1}$ instead of $n$),
\[  \E[
e^{-\lambda L_{ \tau_m}^1} 1_{ \tau_m \leq  r(m)}
     \mid {\cal F}_{n+\tau_{m+1}} ]  \geq  (1-h_m) 
              \E[e^{-\lambda L_{ \tau_m}^1} \mid
       {\cal F}_{n+\tau_{m+1}} ]  \: 1_{\tau_{m+1} \leq  r(m+1)} , \]
where
$ r(m) = \sum_{l=m}^\infty l^2 2^{-l/20}$. 
If we let $m$ be the smallest integer such that $r(m) < 1/4$ and
$h_l < 1$ for all $l \ge m$, then we get 
(\ref{dec14.3}) with $\tau = n + ( \tau_m \wedge
 1/4)$, $\delta_0 = 2^{-m}$ and 
$ c'' = \prod_{l=m}^\infty (1 - h_l) $. 
\qed

\medbreak

If $1 \leq n \leq  m$, let
\[  R_{n,m} = e^{(m-n) \xi}
\E[e^{-\lambda L_m}  
    \mid {\cal F}_n ] 
\hbox { and }  R_{n,m}^* =  e^{(m-n)\xi}
   \E[e^{-\lambda L_m}  1_{G_m}
    \mid {\cal F}_n ] . \]
Then, it follows from
the lemma that there exists constants $c_5,c_6$ such that for
all $m \geq n+1$,
\[       R_{n,m}^* \leq R_{n,m} \leq c_6 R_{n,m}^* , \]
\begin {equation}
\label {eqi}
   c_5 R_{n,n+1}^*  \leq R_{n,m} \leq  c_6 R_{n,n+1}^*
.\end {equation}
This result is used in \cite {LSWan}.

\medbreak

In \cite{Lstrict,LSWan} it is in fact  shown that the limit
$ R_n =  \lim_{m \rightarrow \infty } R_{n,m} $ 
exists and that 
\[    R_{n,m} = R_n [1 +  \epsilon_{n,m}] , \]
where $|\epsilon_{n,m}| \leq c_1 e^{-m c_2}$ and $c_1,c_2$
depend only on $\lambda_0$.
Also, the limit
\[             r = r(\lambda) =
   \lim_{n \rightarrow \infty}  e^{n\xi(2,\lambda)
          } b_n \]
exists and 
\[            b_n = r e^{-n\xi(2,\lambda) } [1 + \epsilon_n] \]
where $|\epsilon_n | \leq c_1 e^{-m c_2}$.

\section {Other exponents and exact values}
\label {final}

The exponents $\xi(2,\lambda)$ comprise just one family
of  Brownian intersection exponents.  The proofs
 apply with minor modifications to these
other exponents.  We review the results here.

Let $\bar p = (p_1,\ldots,p_l)$ be an $l$-tuple of
positive integers and let $\bar \lambda = (\lambda_1,
\ldots,\lambda_l)$ be an $l$-tuple of positive real
numbers.  Let
\[ \BM^{j,k}_t, \;\;\;\;  j=1,2,\ldots,l, \;\;\;\;
    k=1,2,\ldots,p_j , \]
be independent Brownian
motions starting uniformly on $\C_0$.  As
before, let
\[  T^{j,k}_n = \inf\{t > 0 \st  \BM_t^{j,k} \in \C_n \} . \]
For any $j=1, \ldots, l$, define 
$$
{\cal P}_n^j
= \bigcup_{k=1}^{p_j}
\BM^{j,k} [ 0, T_n^{j,l} ].
$$
Let $E_{n,\bar{p}}$ be the event that
the $l$ packets of Brownian motions
  ${\cal P}_n^1, \ldots, {\cal P}_n^l$ 
are disjoint and are ordered clockwise around the origin 
(i.e., %%O-ch
 that their intersection with $\C_n$ are ordered 
clockwise on $\C_n$).
For each $k=1, \ldots, l$, let 
$Z^k_n = Z_n^k ({\cal P}_n^1, \ldots, {\cal P}_n^l)$
denote the probability that a Brownian  
motion $Y$ started uniformly on the unit circle reaches 
$\C_n$ without intersecting $\cup_{j=1}^l {\cal P}_n^j$, 
and in such a way that 
the endpoint of $Y$, $\C_n \cap {\cal P}_n^k$ and $\C_n
\cap {\cal P}_n^{k-1}$  
are ordered 
clockwise on the ${\cal C}_n$ (where ${\cal P}_n^0 = {\cal P}_n^l$).
We then define
$$
b_n ( \lambda_1, p_1, \lambda_2, \ldots, p_l)
=
\expect [1_{E_n, \bar p} \prod_{j=1}^l (Z_n^{j})^{\lambda_j} ].
$$

\begin{theorem}
For every finite integers $M $ and $l$,
there exist constants $0 <
c_1 < c_2 < \infty$ such that the following is true.
For all positive integers $p_1, \ldots , p_l$ that are
smaller than $M$, for all positive reals 
$\lambda_1, \ldots, \lambda_l$ that are smaller than $M$,
 there exists
 $\xi = \xi(\lambda_1, p_1, \ldots,\lambda_l, p_l)$ such
that for all $n \geq 1$,
\[ c_1 e^{-\xi n} \leq 
b_n (\lambda_1, p_1, \ldots, p_l) \leq c_2
   e^{-\xi n} . \]
\end{theorem}

Note (see  \cite {LW1}) that  
$\xi(\lambda_1, p_1, \lambda_2, \ldots, \lambda_l, p_l) 
$ is unchanged if we change the order of the $p$'s and
the $\lambda$'s. Hence, all $b_n$'s (for different
orderings of the $p$'s and the $\lambda$'s) are multiplicative 
constants away from each other.

There are also other exponents called the half-space exponents 
(see \cite {LW1} for a precise definition). The methods of the 
present paper apply also for these exponents. We leave the
detailed 
statement to the interested reader.

Nowhere in this paper have we used the exact values of
the exponents.  Rigorous determination of these values is
the subject of the papers \cite{LSW1,LSW2,LSW3,LSWan}.  In
those papers we prove that
\[  \xi(p_1,\lambda_1,\ldots,p_l,\lambda_l) =
     V[ U (p_1) +U(\lambda_1) +  \cdots +
   U(\lambda_{l})] , \]
where
\[  U(x) = \frac{\sqrt{24 x + 1} -1}{\sqrt{24}} 
\hbox { and }
  V(x) = \frac{6x^2-1}{12}. \]
In particular,
\[ \xi(2,\lambda) = \frac{\lambda}{2} + \frac{11}{24}
   + \frac{5}{24} \sqrt{24 \lambda + 1} . \]

\filbreak
\begingroup
\small 
\parindent=0pt 

\ifhyper\def\email#1{\href{mailto:#1}{\texttt{#1}}}\else
\def\email#1{\texttt{#1}}\fi
\vtop{
\hsize=2.3in
Greg Lawler\\
Department of Mathematics\\
Box 90320\\
Duke University\\
Durham NC 27708-0320, USA\\
\email{jose@math.duke.edu}
}
\bigskip
\vtop{
\hsize=2.3in
Oded Schramm\\
Microsoft Corporation,\\
One Microsoft Way,\\
Redmond, WA 98052; USA\\
\email{schramm@microsoft.com}
}
\bigskip
\vtop{
\hsize=2.3in
Wendelin Werner\\
D\'epartement de Math\'ematiques\\
B\^at. 425\\
Universit\'e Paris-Sud\\
91405 ORSAY cedex, France\\
\email{wendelin.werner@math.u-psud.fr}
}
\endgroup

\filbreak

\end {document}